\documentclass[12pt]{article}

\usepackage{amssymb}
\usepackage{amsmath}
\usepackage{cite}
\usepackage[arrow, matrix, curve]{xy}
\usepackage{tikz}

\begin{document}

\renewcommand{\citeleft}{{\rm [}}
\renewcommand{\citeright}{{\rm ]}}
\renewcommand{\citepunct}{{\rm,\ }}
\renewcommand{\citemid}{{\rm,\ }}

\newcounter{abschnitt}
\newtheorem{satz}{Theorem}
\newtheorem{theorem}{Theorem}[abschnitt]
\newtheorem{koro}[theorem]{Corollary}
\newtheorem{prop}[theorem]{Proposition}
\newtheorem{lem}[theorem]{Lemma}
\newtheorem{conj}[theorem]{Conjecture}
\newtheorem{probl}[theorem]{Problem}
\newtheorem{example}[theorem]{Example}

\newcounter{saveeqn}
\newcommand{\alpheqn}{\setcounter{saveeqn}{\value{abschnitt}}
\renewcommand{\theequation}{\mbox{\arabic{saveeqn}.\arabic{equation}}}}
\newcommand{\reseteqn}{\setcounter{equation}{0}
\renewcommand{\theequation}{\arabic{equation}}}

\hyphenation{convex} \hyphenation{bodies}

\sloppy

\hyphenpenalty=1000

\phantom{a}

\vspace{-2.2cm}

\begin{center}
 \begin{Large} {Binary Operations in Spherical Convex Geometry} \\[0.7cm] \end{Large}

\begin{large} Florian Besau and Franz E.
Schuster \end{large}
\end{center}

\vspace{-1.1cm}

\begin{quote}
\footnotesize{ \vskip 1truecm\noindent {\bf Abstract.}
Characterizations of binary operations between convex bodies on
the Euclidean unit sphere are established. The main result shows
that the convex hull is essentially the only non-trivial
projection covariant operation between pairs of convex bodies
contained in open hemispheres. Moreover, it is proved that any
continuous and projection covariant binary operation between all
proper spherical convex bodies must be trivial. }
\end{quote}

\vspace{0.6cm}

\centerline{\large{\bf{ \setcounter{abschnitt}{1}
\arabic{abschnitt}. Introduction}}} \alpheqn

\vspace{0.6cm}

In recent years it has been explained why a number of fundamental
notions from convex geometric analysis really do have a special
place in the theory. For example, Blaschke's classical affine and
centro-affine surface areas were given characterizations by Ludwig
and Reitzner \textbf{\cite{Ludwig:2010}} and Haberl and Parapatits
\textbf{\cite{HabPap2013}} as unique valuations satisfying
certain invariance properties; polar duality and the Legendre
transform were characterized by B\"or\"oczky and Schneider
\textbf{\cite{BoeroecSchn2008}} and Artstein-Avidan and Milman
\textbf{\cite{ArtsteinMilman2009}}, respectively. These and other
results of the same nature (see also, e.g.,
\textbf{\cite{Haberl:2006, Ludwig:2003, Ludwig:2005, Schu2010,
SchWan2012, Wan2011 }}) not only show that the notions under
consideration are characterized by a surprisingly small number of
basic properties but also led to the discovery of seminal new
notions.

Gardner, Hug, and Weil \textbf{\cite{GardHugWeil2013}} initiated
a new line of research whose goal is to enhance our understanding
of the fundamental characteristics of known {\it binary}
operations between sets in Euclidean geometry (see also
\textbf{\cite{GardParapSchu2014}}). Their main focus is on
operations which are projection covariant, that is, the \linebreak
operation can take place before or after projection onto linear
subspaces, with the same effect. One impressive example of the
results obtained in \textbf{\cite{GardHugWeil2013}} is a
characterization of the classical Minkowski addition between
convex bodies (compact convex sets) in $\mathbb{R}^n$ as the only
projection covariant operation which also satisfies the identity
property. In fact, a characterization of all projection covariant
operations between origin-symmetric convex bodies was established
in \textbf{\cite{GardHugWeil2013}}, by proving that such
operations are precisely those given by so-called $M$-addition
(see Section 3 for precise definitions). This little-known
addition was later shown in \textbf{\cite{GardHugWeil2014}} to be
intimately related to Orlicz addition, a recent important
generalization of Minkowski addition.

\pagebreak

The Brunn--Minkowski theory, which arises from combining volume
{\it and} Minkowski addition, lies at the very core of classical
Euclidean convexity and provides a unifying framework for various
extremal and uniqueness problems for convex bodies in
$\mathbb{R}^n$ (see, e.g., \textbf{\cite{Gardner:2006,
gruber2009, Schneider:1993}}). In contrast, the geometry of
spherical convex sets is much less well understood. Although
certain aspects, like the integral geometry of spherical convex
sets (see \textbf{\cite{alesker2007, alesker2010,
bernigfusol2014, Glasauer1996, Santalo:2004, Klain:1997}}), have witnessed
considerable progress, contributions to spherical convexity are
rather scattered (see \textbf{\cite{Aubrun:2004, Burago:1988,
Figiel:1977, Gao:2003, Gardner:2002, Gerhardt:2007, Oliker:2007,
Santalo:1950, Santalo:1980, schneiderweil, Yaskin2006}}). The
reason for this might be that so far no natural analogue of
Minkowski addition is available on the sphere. (For an attempt to
remedy this see \textbf{\cite{leichtweiss2012}}.)

In this article we start a systematic investigation of binary
operations between convex bodies (that is, closed convex sets) on
the Euclidean unit sphere with a focus on operations which are
covariant under projections onto great subspheres. We prove that
all continuous such operations between proper spherical convex
bodies are trivial. More importantly, our main result shows that
the convex hull is essentially the only non-trivial projection
covariant operation between pairs of convex bodies contained in
open hemispheres. The picture changes drastically when operations
between convex bodies in a {\it fixed} open hemisphere are
considered. In this case, we establish a one-to-one
correspondence between binary operations on spherical convex
bodies that are projection covariant with respect to the center
of the hemisphere, and projection covariant operations on convex
bodies in $\mathbb{R}^n$.

\vspace{1cm}

\centerline{\large{\bf{ \setcounter{abschnitt}{2}
\arabic{abschnitt}. Statement of principal results}}}

\reseteqn \alpheqn

\vspace{0.6cm}

\newcommand{\R}{\mathbb{R}}
\renewcommand{\S}{\mathbb{S}}
\newcommand{\K}{\mathcal{K}}
\renewcommand{\O}{u}
\newcommand{\Kps}{\K_{\O}(\S^n)}
\newcommand{\Kp}{\K_{\circ}(\S^n)}

Let $\mathbb{S}^n$ denote the $n$-dimensional Euclidean unit
sphere. Throughout the article we assume that $n \geq 2$. The
usual \emph{spherical distance} between points on $\mathbb{S}^n$
is given by $d(u,v) = \arccos(u \cdot v)$, $u, v \in
\mathbb{S}^n$. For $\lambda > 0$ and $A \subseteq \mathbb{S}^n$,
we write $A_\lambda$ for the set of all points with distance at
most $\lambda$ from $A$. The \emph{Hausdorff distance} between
closed sets $A, B \subseteq \mathbb{S}^n$ is then given by
\[\delta_s(A,B) = \min\left\{0\leq \lambda \leq \pi: A\subseteq B_\lambda \text{ and } B\subseteq A_\lambda\right\}.\]
A set $A \subseteq \mathbb{S}^n$ is called \emph{(spherical)
convex} if
\[\mathrm{rad}\,A =\left\{\lambda x: \lambda \geq 0, x \in A
\right\}\subseteq\R^{n+1}\] is convex. We say $K \subseteq
\mathbb{S}^n$ is a \emph{convex body} if $K$ is closed and convex.
Let $\K(\mathbb{S}^n)$ denote the space of convex bodies in
$\mathbb{S}^n$ with the Hausdorff distance.

\pagebreak

We call $K\in\K(\mathbb{S}^n)$ a \emph{proper} convex body if $K$
is contained in an open hemisphere and we write
$\mathcal{K}^p(\mathbb{S}^n)$ for the subspace of
$\K(\mathbb{S}^n)$ of all proper convex bodies. For fixed $u \in
\mathbb{S}^n$ we denote by $\mathcal{K}^p_u(\mathbb{S}^n)$ the
subspace of (proper) convex bodies that are contained in the open
hemisphere centered at $u$. Then
\[ \mathcal{K}^p(\mathbb{S}^n) = \bigcup_{\substack{u \in\S^n}}\mathcal{K}_u^p(\mathbb{S}^n).\]

The \emph{convex hull} of $A \subseteq \mathbb{S}^n$ is the
intersection of all convex sets in $\mathbb{S}^n$ that contain
$A$. Note that, for $K, L \in \mathcal{K}(\mathbb{S}^n)$, we have
$\mathrm{conv}(K \cup L) \in \mathcal{K}(\mathbb{S}^n).$

For $0 \leq k \leq n$, a \emph{$k$-sphere} $S$ is a
$k$-dimensional great sub-sphere of $\mathbb{S}^n$, that is, the
intersection of a $(k+1)$-dimensional linear subspace $V
\subseteq \R^{n+1}$ with $\S^n$. Clearly, every $k$-sphere $S$ is
convex. For $K \in \mathcal{K}(\mathbb{S}^n)$, the \emph{spherical
projection} $K|S$ is defined by
\[K|S = \mathrm{conv}\left(K\cup S^\circ\right)\cap S = \left(\mathrm{rad}(K)|V\right)\cap
\S^n,\] where $S=V\cap \S^n$ and $S^\circ$ is the $(n-k-1)$-sphere
orthogonal to $S$, that is, $S^\circ=V^\bot\cap \S^n$.

For fixed $u \in \S^n$ we call a binary operation $*\colon
\mathcal{K}^p(\mathbb{S}^n) \times \mathcal{K}^p(\mathbb{S}^n)
\rightarrow \mathcal{K}^p(\mathbb{S}^n)$ \emph{$\O$-projection
covariant} if for all $k$-spheres $S$, $0 \leq k \leq n - 1$, with
$u \in S$ and for all $K$, $L \in \mathcal{K}_u^p(\mathbb{S}^n)$,
we have
\[(K|S) * (L|S) = (K*L)|S.\]
We call $*$ \emph{projection covariant} if $*$ is $\O$-projection
covariant for all $\O\in\S^n$.

\vspace{0.2cm}

The main objective of this article is to characterize projection
covariant operations between spherical convex bodies. Our first
result shows that such operations between \emph{all} proper
convex bodies in $\mathbb{S}^n$ are of a very simple form.

\begin{satz} \label{thm:1} An operation $*\colon \mathcal{K}^p(\mathbb{S}^n) \times \mathcal{K}^p(\mathbb{S}^n) \to \mathcal{K}^p(\mathbb{S}^n)$
between proper convex bodies is projection covariant and
continuous with respect to the Hausdorff metric if and only if
either $K * L = K$, or $K * L = -K$, or $K * L = L$, or $K * L =
-L$ for all $K, L \in \mathcal{K}^p(\mathbb{S}^n)$.
\end{satz}

We call the binary operations from Theorem \ref{thm:1}
\emph{trivial}. As the following example shows, the continuity
assumption in Theorem \ref{thm:1} cannot be omitted.

\pagebreak

\vspace{0.3cm}

\noindent {\bf Example:}

\vspace{0.1cm}

\noindent Consider the set $\mathcal{C} \subset
\mathcal{K}^p(\mathbb{S}^n) \times \mathcal{K}^p(\mathbb{S}^n)$
of all pairs $(K,L)$ such that both $K$ and $L$ are contained in
some open hemisphere, that is,
\begin{align*}
    \mathcal{C}=\bigcup_{\substack{\O \in \S^n}} \left( \mathcal{K}^p_u(\mathbb{S}^n) \times \mathcal{K}^p_u(\mathbb{S}^n)  \right).
\end{align*}
Define an operation $*\colon \mathcal{K}^p(\mathbb{S}^n) \times
\mathcal{K}^p(\mathbb{S}^n) \to \mathcal{K}^p(\mathbb{S}^n)$ by
\begin{align*}
    K*L =\begin{cases}
        K &\text{if $(K,L)\in\mathcal{C}$},\\
        L &\text{if $(K,L)\notin\mathcal{C}$.}
    \end{cases}
\end{align*}
Clearly, $*$ is not continuous but by our definition it is
projection covariant.

\vspace{0.3cm}

The proof of Theorem \ref{thm:1} relies on ideas of Gardner, Hug,
and Weil. The critical tool to transfer their techniques to the
sphere is the gnomonic projection (see Section 4) which
establishes the following correspondence \linebreak between
projection covariant operations on $\mathcal{K}(\mathbb{R}^n)$,
the space of compact, convex sets in $\mathbb{R}^n$, and
$u$-projection covariant operations on
$\mathcal{K}_u^p(\mathbb{S}^n)$:

\begin{satz}\label{thm:2}
For every fixed $u \in \mathbb{S}^n$, there is a one-to-one
correspondence between $u$-projection covariant operations
$*\colon \mathcal{K}_u^p(\mathbb{S}^n) \times
\mathcal{K}_u^p(\mathbb{S}^n) \to \mathcal{K}_u^p(\mathbb{S}^n)$
and projection covariant operations $\overline{*}\colon
\K(\R^n)\times\K(\R^n)\to\K(\R^n)$. Moreover, every such
$\O$-projection covariant operation $*$ is continuous in the
Hausdorff metric.
\end{satz}

Note that by Theorem \ref{thm:2} every projection covariant
operation $*$ on $\mathcal{C}$ is also automatically continuous.

\vspace{0.1cm}

Finally, as our main result, we prove that the only {\it
non-trivial} projection covariant operation on the set
$\mathcal{C}$ is essentially the spherical convex hull.

\begin{satz}\label{thm:3}
An operation $*\colon \mathcal{C} \to
\mathcal{K}^p(\mathbb{S}^n)$ is non-trivial and projection
covariant if and only if either $K * L = \mathrm{conv}(K\cup L)$
or $K * L = - \mathrm{conv}(K \cup L)$ for all $(K,L) \in
\mathcal{C}$.
\end{satz}

After briefly recalling the background material on convex bodies
in $\mathbb{R}^n$ in Section 3, we discuss the geometry of spherical convex sets in Section 4 and use the gnomonic
projection to prove Theorem \ref{thm:2}.
Sections 5 and 6 contain the proofs of Theorems \ref{thm:1} and \ref{thm:3}. Motivated by investigations of Gardner, Hug, and Weil
\textbf{\cite{GardHugWeil2013}} in $\mathbb{R}^n$, we discuss
{\it section} covariant operations between spherical {\it star
sets} in the concluding section of the article.

\pagebreak

\centerline{\large{\bf{ \setcounter{abschnitt}{3}
\arabic{abschnitt}. Background material from Euclidean
convexity}}}

\reseteqn \alpheqn

\vspace{0.6cm}

In this section we collect basic material about convex bodies in
$\mathbb{R}^n$. As a general reference for these facts we
recommend \textbf{\cite{Schneider:1993}}. We also recall the
definition of the $L_p$ Minkowski addition and, more generally,
the $M$-addition of convex bodies as well as their characterizing
properties established in \textbf{\cite{GardHugWeil2013}}.

The standard orthonormal basis for $\mathbb{R}^n$ will be $\{e_1,
\ldots, e_n\}$. Otherwise, we usually denote the coordinates of
$x \in \mathbb{R}^n$ by $x_1, \ldots, x_n$. We write $B^n$ for
the Euclidean unit ball in $\mathbb{R}^n$. We call a subset of
$\mathbb{R}^n$ \emph{$1$-unconditional} if it is symmetric with
respect to each coordinate hyperplane.

Let $\mathcal{K}_e(\mathbb{R}^n)$ be the set of origin symmetric
convex bodies and let $\mathcal{K}_{\mathrm{o}}(\mathbb{R}^n)$
denote the set of convex bodies containing the origin.

A convex body $K \in \mathcal{K}(\mathbb{R}^n)$ is uniquely
determined by its \emph{support function} defined by
\[h(K,x) = \max\{x \cdot y: y \in K\}, \qquad x \in \mathbb{R}^n.\]
We will sometimes also use $h_K$ to denote the support function
of $K \in \mathcal{K}(\mathbb{R}^n)$. Support functions are
\emph{$1$-homogeneous}, that is, $h(K,\lambda x) = \lambda
h(K,x)$ for all $x \in \mathbb{R}^n$ and $\lambda > 0$, and are
therefore often regarded as functions on $\mathbb{S}^{n-1}$. They
are also \emph{subadditive}, that is,  $h(K,x + y) \leq h(K,x) +
h(K,y)$ for all $x, y \in \mathbb{R}^n$. Conversely, every
1-homogeneous and subadditive function on $\mathbb{R}^n$ is the
support function of a convex body. Clearly, $K \in
\mathcal{K}_e(\mathbb{R}^n)$ if and only if $h(K,\cdot)$ is even.

The \emph{Minkowski sum} of subsets $X$ and $Y$ of $\mathbb{R}^n$
is defined by
\[X + Y = \{x + y : x \in X, y \in Y \}.\]
If $K, L \in \mathcal{K}(\mathbb{R}^n)$, then $K + L$ can be
equivalently defined as the convex body such that
\[h(K+L,\cdot) = h(K,\cdot) + h(L,\cdot).\]

The \emph{Hausdorff distance} $\delta(X,Y)$ between compact
subsets $X$ and $Y$ of $\mathbb{R}^n$ is defined by
\[\delta(X,Y) = \min\{\lambda \geq 0: X \subseteq Y + \lambda\,B^n \mbox{ and } Y \subseteq X + \lambda\,B^n \}.  \]
If $K, L \in \mathcal{K}(\mathbb{R}^n)$, then $\delta(K,L)$ can be
alternatively defined by
\begin{equation} \label{hausdsupp}
\delta(K,L) = \|h(K,\cdot) - h(L,\cdot)\|_{\infty},
\end{equation}
where $\|\cdot\|_{\infty}$ denotes the $L_{\infty}$ norm on
$\mathbb{S}^{n-1}$.

\pagebreak

For $1 < p \leq \infty$, the \emph{$L_p$ Minkowski sum} of convex
bodies $K, L \in \mathcal{K}_{\mathrm{o}}(\mathbb{R}^n)$ was first
defined by Firey \textbf{\cite{Firey1962}} by
\[h(K +_p L,\cdot)^p = h(K,\cdot)^p + h(L,\cdot)^p,  \]
for $p < \infty$, and by
\[h(K +_{\infty} L,\cdot) = \max\{h(K,\cdot),h(L,\cdot)\}.  \]
Note that $K +_{\infty} L$ is just the usual convex hull in
$\mathbb{R}^n$ of $K$ and $L$.

Lutwak \textbf{\cite{Lutwak:1993, Lutwak:1996}} showed that the
$L_p$ Minkowski addition leads to a very powerful extension
of the classical Brunn--Minkowski theory. Since the 1990's this
$L_p$ Brunn--Minkowski theory has provided new tools for attacks
on major unsolved problems and consolidated connections between
convex geometry and other fields (see, e.g.,
\textbf{\cite{BLYZ:2013, Ludwig:2003, Lutwak:2000, Lutwak:2000a,
LYZ:2004, Ryabogin2004, stancu2002, Web2013, Werner2012,
Werner2008}} and the references therein). An extension of the
$L_p$ Minkowski addition to arbitrary sets in $\mathbb{R}^n$ was
given only recently in \textbf{\cite{LYZ2012}}.

An even more general way of combining two subsets of
$\mathbb{R}^n$ is the still more recent \emph{$M$-addition}: If
$M$ is an arbitrary subset of $\mathbb{R}^2$, then the
\emph{$M$-sum} of $X, Y \subseteq \mathbb{R}^n$ is defined by
\begin{equation} \label{defmadd}
X \oplus_M Y = \bigcup_{(a,b) \in M}\! a\,X + b\,Y = \left \{ax +
by: (a,b) \in M, x \in X, y \in Y \right \}.
\end{equation}

Protasov \textbf{\cite{protasov1997}} first introduced
$M$-addition for centrally symmetric convex bodies and a
$1$-unconditional convex body $M$ in $\mathbb{R}^2$. He also
proved that $\oplus_M: \mathcal{K}_e(\mathbb{R}^n) \times
\mathcal{K}_e(\mathbb{R}^n) \rightarrow
\mathcal{K}_e(\mathbb{R}^n)$ for such $M$.

Gardner, Hug and Weil \textbf{\cite{GardHugWeil2013}}
rediscovered $M$-addition in the more general form (\ref{defmadd})
in their investigation of projection covariant binary operations
between convex bodies in $\mathbb{R}^n$. Among several results on
this seminal operation, they proved the following:

\begin{theorem} {\bf ( \hspace{-0.3cm} \cite{GardHugWeil2013})} Let $M \subseteq \mathbb{R}^2$.
Then $\oplus_M: \mathcal{K}(\mathbb{R}^n) \times
\mathcal{K}(\mathbb{R}^n) \rightarrow \mathcal{K}(\mathbb{R}^n)$
if and only if $M \in \mathcal{K}(\mathbb{R}^2)$ and $M$ is
contained in one of the 4 quadrants of $\mathbb{R}^2$. In this
case, let $\varepsilon_i = \pm 1$, $i=1,2$, denote the sign of
the $i$th coordinate of a point in the interior of this quadrant
and let
\[M^+=\{(\varepsilon_1a,\varepsilon_2b): (a,b) \in M\}  \]
be the reflection of $M$ contained in $[0,\infty)^2$. If $K, L
\in \mathcal{K}(\mathbb{R}^n)$, then
\begin{equation} \label{suppmadd}
h_{K \oplus_M L}( x) = h_{M^+}(h_{\varepsilon_1K}(x),h_{\varepsilon_2L}(x)), \qquad x \in
\mathbb{R}^n.
\end{equation}
\end{theorem}

\pagebreak

\noindent {\bf Example:}

\vspace{0.1cm}

\noindent For some $1 \leq p \leq \infty$, let
\[M = \{(a,b) \in [0,1]^2: a^{p'} + b^{p'} \leq 1\},   \]
where $1/p + 1/p' = 1$. Then $\oplus_M = +_p$ is $L_p$ Minkowski
addition on $\mathcal{K}_{\mathrm{o}}(\mathbb{R}^n)$.

\vspace{0.4cm}

The following basic properties of $M$-addition are of particular
interest for us. They are immediate consequences of either
definition (\ref{defmadd}) or (\ref{suppmadd}).

\begin{prop} Suppose that $M \in \mathcal{K}(\mathbb{R}^2)$ is contained in $[0,\infty)^2$. Then $\oplus_M: \mathcal{K}(\mathbb{R}^n) \times
\mathcal{K}(\mathbb{R}^n) \rightarrow \mathcal{K}(\mathbb{R}^n)$
has the following properties:
\begin{itemize}
\item \emph{Continuity}\\[0.1cm]
$K_i \to K$, $L_i \to L$ implies $K_i \oplus_M L_i \to K \oplus_M
L$ as $i \rightarrow \infty$ in the Hausdorff metric;
\item \emph{$\mathrm{GL}(n)$ covariance} \\[0.1cm]
$(AK) \oplus_M (AL) = A(K \oplus_M L)$ for all $A \in \mathrm{GL}(n)$;
\item \emph{Projection covariance}\\[0.1cm] $(K|V) \oplus_M (L|V) = (K \oplus_M L)|V$ for every linear subspace $V$ of $\mathbb{R}^n$.
\end{itemize}
\end{prop}

It is easy to show that continuity and $\mathrm{GL}(n)$
covariance imply projection covariance. That the converse
statement is also true, follows from a deep result of Gardner,
Hug, and Weil which states the following:

\begin{theorem} \label{thm:gardner} {\bf ( \hspace{-0.2cm} \cite{GardHugWeil2013})}
An operation $*: \mathcal{K}(\mathbb{R}^n) \times
\mathcal{K}(\mathbb{R}^n) \to \mathcal{K}(\mathbb{R}^n)$ is
projection covariant if and only if there exists a nonempty closed
convex set $\overline{M}$  in $\mathbb{R}^4$ such that, for all $K$, $L \in \mathcal{K}(\mathbb{R}^n)$,
\begin{equation} \label{suppmadd4}
h_{K * L}(x) =
h_{\overline{M}}(h_{-K}(x),h_K(x),h_{-L}(x),h_L(x)), \qquad x \in
\mathbb{R}^n.
\end{equation}
Consequently, every such operation is continuous and
$\mathrm{GL}(n)$ covariant.
\end{theorem}

Note that it is an open problem whether the binary operation on
$\mathcal{K}(\mathbb{R}^n)$ defined by (\ref{suppmadd4}) is
$M$-addition for some (convex) subset $M$ of $\mathbb{R}^2$.
However, Gardner, Hug, and Weil \textbf{\cite{GardHugWeil2013}}
proved that an operation between \emph{$o$-symmetric} convex
bodies is projection covariant if and only if it is $M$-addition
for some $1$-unconditional convex body in $\mathbb{R}^2$.

\pagebreak

\centerline{\large{\bf{ \setcounter{abschnitt}{4}
\arabic{abschnitt}. The gnomonic projection}}}

\reseteqn \alpheqn \setcounter{theorem}{0}

\vspace{0.6cm}

In the following we discuss basic facts about spherical convex
sets. In particular, we recall the definition of spherical
support functions of proper convex bodies in $\mathbb{S}^n$. The
second part of this section is devoted to the gnomonic projection.
After establishing the basic properties of this critical tool, we conclude this section with the proof of Theorem \ref{thm:2}.

For the following alternative definitions of proper convex bodies
in $\mathbb{S}^n$, we refer to \textbf{\cite{danzeretal1963}}.

\begin{prop} The following statements about $K \subseteq \mathbb{S}^n$ are equivalent:
\begin{enumerate}
\item[(a)] The set $K$ is a proper convex body.
\item[(b)] The set $K$ is an intersection of open hemispheres.
\item[(c)] There are no antipodal points in $K$ and for every two points $u,v \in K$,
the minimal geodesic connecting $u$ and $v$ is contained in $K$.
\end{enumerate}
\end{prop}

Although we will make no use of this fact, we remark, that a set
$K \subseteq \mathbb{S}^n$ is a convex body if and only if $K$ is
the intersection of \emph{closed} hemispheres.

\vspace{0.2cm}

Next we introduce spherical support functions of proper convex
sets contained in a fixed hemisphere (cf.\
\textbf{\cite{leichtweiss2012}} for a related construction). To
this end, for $u \in \mathbb{S}^n$, let $\mathbb{S}_u^+$ denote
the open hemisphere with center in $u$ and let $\mathbb{S}_u$ be
the boundary of $\mathbb{S}_u^+$, that is,
\[\mathbb{S}_u^+ = \{v \in \mathbb{S}^n: u \cdot v > 0\} \qquad \mbox{and} \qquad \mathbb{S}_u = \{v \in \mathbb{S}^n: u \cdot v = 0\}.  \]
For non-antipodal $u, v \in \mathbb{S}^n$, we write $S_{u,v}$ for
the unique great circle containing $u$ and $v$.

\vspace{0.3cm}

\noindent {\bf Definition} For $u \in \mathbb{S}^n$ and a proper convex body $K \in
\mathcal{K}_u^p(\mathbb{S}^n)$, the \emph{spherical support
function} $h_u(K,\cdot): \mathbb{S}_u \rightarrow \left ( -
\frac{\pi}{2},\frac{\pi}{2} \right )$ of $K$ is defined by
\[h_u(K,v) = \max\{\mathrm{sgn}(v \cdot w)\,d(u,w|S_{u,v}): w \in K\}.  \]

\vspace{0.2cm}

Recall that the (Euclidean) support function of a convex body $L$ in
$\mathbb{R}^n$ encodes the signed distances of the supporting
planes to $L$ from the origin. In other words, we have for every $v
\in \mathbb{S}^{n-1}$,
\[L|\mathrm{span}\{v\} = \{tv: t \in [-h(K,-v),h(K,v)]\}.\]
The intuitive meaning of the spherical support function of a
proper convex body $K \in \mathcal{K}_u^p(\mathbb{S}^n)$ is
similar. It yields the oriented angle between $u$ and the supporting
$(n-1)$-spheres to $K$. More precisely, we have for every $v \in
\mathbb{S}_u$,
\[K|S_{u,v} = \{u \cos \alpha + v \sin \alpha: \alpha \in [-h_u(K,-v),h_u(K,v)]\}. \]
In particular, for $K, L \in \mathcal{K}_u^p(\mathbb{S}^n)$, $K
\subseteq L$ if and only if $h_u(K,\cdot) \leq h_u(L,\cdot)$.

In the following we denote by $\mathbb{R}^n_u$ (instead of
$u^{\bot}$) the hyperplane in $\mathbb{R}^{n+1}$ orthogonal to $u
\in \mathbb{S}^n$.

\vspace{0.3cm}

\noindent {\bf Definition} For $u \in \mathbb{S}^n$, the
\emph{gnomonic projection} $g_u: \mathbb{S}_u^+ \rightarrow
\mathbb{R}^n_u$ is defined by
\[g_u(v) = \frac{v}{u \cdot v} - u.  \]

\vspace{0.2cm}

In the literature, the gnomonic projection is often considered as a map
to the tangent plane at $u$. However, for our purposes it is more convenient if the range of $g_u$ contains the origin.

In the following lemma we collect a number of well-known properties of the gnomonic projection which are
immediate consequences of its definition.

\begin{lem} \label{gnombasic1} For $u \in \mathbb{S}^n$, the following statements hold:
\begin{enumerate}
\item[(a)] The gnomonic projection $g_u: \mathbb{S}_u^+ \rightarrow
\mathbb{R}^n_u$ is a bijection with inverse
\[g_u^{-1}(x) = \frac{x + u}{\|x + u\|}, \qquad x \in \mathbb{R}^n_u.\]
\item[(b)] If $S \subseteq \mathbb{S}^n$ is a $k$-sphere, $0 \leq k \leq n - 1$, such that $S \cap \mathbb{S}_u^+$ is non-empty, then $g_u(S \cap \mathbb{S}_u^+)$ is a
$k$-dimensional affine subspace of $\mathbb{R}^n_u$. Conversely,
$g_u^{-1}$ maps $k$-dimensional affine subspaces of
$\mathbb{R}^n_u$ to $k$-hemispheres in $\mathbb{S}_u^+$.
\item[(c)] The gnomonic projection maps $\mathcal{K}_u^p(\mathbb{S}^n)$
bijectively to $\mathcal{K}(\mathbb{R}_u^n)$.
\end{enumerate}
\end{lem}

If $K \in \mathcal{K}_u^p(\mathbb{S}^n)$, then, by Lemma
\ref{gnombasic1} (c), the set $g_u(K)$ is a convex body in
$\mathcal{K}(\mathbb{R}_u^n)$. The next lemma relates the
(Euclidean) support function of $g_u(K)$ with the spherical
support function of $K$.

\begin{lem} \label{suppsupp}
For $u \in \mathbb{S}^n$ and every $K \in
\mathcal{K}_u^p(\mathbb{S}^n)$, we have
\[h(g_u(K),v) = \tan h_u(K,v), \qquad v \in \mathbb{S}_u.  \]
In particular, $K$ is uniquely determined by $h_u(K,\cdot)$.
\end{lem}
{\it Proof.} For $v \in \mathbb{S}_u$ and $w \in \mathbb{S}_u^+$,
an elementary calculation shows that
\[\frac{v \cdot w}{u \cdot w} = \tan(\mathrm{sgn}(v \cdot w)\,d(u,w|S_{u,v})).  \]
Therefore, the definition of $g_u$ and the monotonicity of the
tangent yield
\[h(g_u(K),v) = \max_{x \in g_u(K)}\{ v \cdot x\} = \max_{w \in K} \left
\{\frac{v \cdot w}{u \cdot w} \right \}  = \tan h_u(K,v).
\]

\vspace{-0.5cm}

\hfill $\blacksquare$

\vspace{0.4cm}

By Lemma \ref{suppsupp}, a function $h: \mathbb{S}_u \rightarrow
\left ( - \frac{\pi}{2},\frac{\pi}{2} \right )$ is the spherical
support function of a convex body $K \in
\mathcal{K}_u^p(\mathbb{S}^n)$ if and only if the $1$-homogeneous
extension of $\tan h$ to $\mathbb{R}_u^n$ is the support function
of a convex body in $\mathbb{R}_u^n$.

\vspace{0.2cm}

Using spherical support functions, we define a metric $\gamma_u$
on $\mathcal{K}_u^p(\mathbb{S}^n)$ by
\[ \gamma_u(K,L) = \max_{v\in \mathbb{S}_u} |h_u(K,v)-h_u(L,v)|.\]
Since for $K \in \mathcal{K}_u^p(\mathbb{S}^n)$ and $\varepsilon >
0$, the set $K_{\varepsilon}$ of all points with distance at most
$\varepsilon$ from $K$ is \emph{not} necessarily convex, it is not
difficult to see that the restriction of $\delta_s$ to
$\mathcal{K}_u^p(\mathbb{S}^n)$ does \emph{not} coincide with
$\gamma_u$ (in contrast to the Euclidean setting). However, our
next result shows that $\gamma_u$ and $\delta_s$ induce the same
topology on $\mathcal{K}_u^p(\mathbb{S}^n)$. Since we could not
find a reference for this basic result, we include a proof for
the readers convenience.

\begin{prop} \label{schweregeb} For $u \in \mathbb{S}^n$, the metrics $\gamma_u$ and $\delta_s$
induce the same topology on $\mathcal{K}_u^p(\mathbb{S}^n)$.
\end{prop}
{\it Proof.} Let $K \in \mathcal{K}_u^p(\mathbb{S}^n)$ and
$\varepsilon > 0$ sufficiently small. We denote by
$B_{\gamma_u}(K,\varepsilon)$ the metric ball with respect to
$\gamma_u$ of radius $\varepsilon$ and center $K$ and
$B_{\delta_s}(K,\varepsilon)$ is defined similarly. We first show
that there exists $r(K,\varepsilon) > 0$ such that
\begin{equation} \label{ball2}
B_{\delta_s}(K,r(K,\varepsilon)) \subseteq
B_{\gamma_u}(K,\varepsilon).
\end{equation}
To this end, let again $w \in \mathbb{S}_u^+$. Since
$\{w\}_{\varepsilon}$ is a spherical cap of radius $\varepsilon$,
it is not difficult to show that
\[\max_{v \in \mathbb{S}_u} (h_u(\{w\}_{\varepsilon},v) - h_u(\{w\},v)) = \arcsin\left ( \frac{\sin \varepsilon}{u\cdot w}  \right ), \]
where this maximum is attained for $v \in \mathbb{S}_u \cap
\mathbb{S}_w$. Therefore, if we define
\[c(w,\varepsilon) = \arcsin (u \cdot w\,\sin \varepsilon  )  ,\]
then
\begin{equation} \label{luke17}
h_u(\{w\}_{c(w,\varepsilon)},v) \leq h_u(\{w\},v) + \varepsilon
\end{equation}
for all $v \in \mathbb{S}_u$. We now define
\[r(K,\varepsilon) = \min_{w \in K} c(w,\varepsilon/2 ).  \]
Note that, by the compactness of $K$, we have $r(K,\varepsilon) >
0$. Since (\ref{luke17}) holds for all $w \in \mathbb{S}_u^+$, we
obtain
\begin{eqnarray*}
\max_{w \in K} h_u(\{w\}_{r(K,\varepsilon)},v) & \leq &  \max_{w \in K} h_u(\{w\}_{c(w,\varepsilon/2)},v) \\
& \leq & \max_{w \in K} h_u(\{w\},v) + \frac{\varepsilon}{2} \leq
h_u(K,v) + \frac{\varepsilon}{2}
\end{eqnarray*}
for all $v \in \mathbb{S}_u$. Using
\[\max_{w \in K} h_u(\{w\}_{r(K,\varepsilon)},v) = h_u(\mathrm{conv}(K_{r(K,\varepsilon)}),v), \qquad v \in \mathbb{S}_u,   \]
we conclude that
\begin{equation} \label{hansolo17}
h_u(\mathrm{conv}(K_{r(K,\varepsilon)}),v) \leq h_u(K,v) +
\frac{\varepsilon}{2}
\end{equation}
for all $v \in \mathbb{S}_u$. Moreover,
\begin{equation} \label{obiwan17}
r(\{w\}_{r(\{w\},\varepsilon)},2\varepsilon) = \min_{w' \in
\{w\}_{c(w,\varepsilon/2)}}c(w',\varepsilon) \geq
c(w,\varepsilon/2).
\end{equation}
This follows from an elementary calculation and the fact, that
$c(w',\varepsilon)$ attains its minimum in
$\{w\}_{c(w,\varepsilon/2)}$ when $d(w',u) =
d(w,u)+c(w,\varepsilon/2)$.

Now, let  $L \in \mathcal{K}_u^p(\mathbb{S}^n)$ such that
$\delta_s(K,L) \leq r(K,\varepsilon)$. Then, from
\begin{equation} \label{leia17}
L \subseteq K_{r(K,\varepsilon)} \subseteq
\mathrm{conv}(K_{r(K,\varepsilon)})
\end{equation}
and (\ref{hansolo17}), we obtain on the one hand
\[h_u(L,v) \leq h_u(\mathrm{conv}(K_{r(K,\varepsilon)}),v) \leq h_u(K,v) +
\frac{\varepsilon}{2} \leq h_u(K,v) + \varepsilon \] for all $v
\in \mathbb{S}_u$. On the other hand, from (\ref{leia17}) and
(\ref{obiwan17}) we deduce that
\[r(L,2\varepsilon) \geq\! \min_{w \in K_{r(K,\varepsilon)}}\! c(w,\varepsilon) \geq \min_{w \in K} r(\{w\}_{r(\{w\},\varepsilon)},2\varepsilon) \geq \min_{w \in K} c(w,\varepsilon/2) = r(K,\varepsilon)  \]
and, thus, $K \subseteq L_{r(K,\varepsilon)} \subseteq
L_{r(L,2\varepsilon)} \subseteq
\mathrm{conv}(L_{r(L,2\varepsilon)})$. Consequently, by
(\ref{hansolo17}),
\[h_u(K,v) \leq h_u(\mathrm{conv}(L_{r(L,2\varepsilon)}),v) \leq h_u(L,v) + \varepsilon \]
for all $v \in \mathbb{S}_u$, which concludes the proof of
(\ref{ball2}).

\vspace{0.2cm}

It remains to show that there also exists
$\overline{r}(K,\varepsilon)
> 0$ such that
\begin{equation} \label{ball1}
B_{\gamma_u}(K,\overline{r}(K,\varepsilon)) \subseteq
B_{\delta_s}(K,\varepsilon).
\end{equation}
To this end, let again $w \in \mathbb{S}_u^+$. By our definition
of spherical support functions, we have for sufficiently small
$\lambda > 0$,
\[\min_{v \in \mathbb{S}_u} (h_u(\{w\}_{\lambda},v) - h_u(\{w\},v)) = \lambda,   \]
where this minimum is attained for $v \in \mathbb{S}_u \cap
S_{u,w}$. Consequently, we obtain
\[h_u(\{w\},v) + \lambda \leq h_u(\{w\}_{\lambda},v)  \]
for all $v \in \mathbb{S}_u$. Since this holds for all $w \in
\mathbb{S}_u^+$, we conclude that
\[h_u(K,v) + \lambda = \max_{w \in K} h_u(\{w\},v) + \lambda \leq  \max_{w \in K} h_u(\{w\}_{\lambda},v) = h_u(\mathrm{conv}(K_{\lambda}),v)     \]
for all $v \in \mathbb{S}_u$. Therefore, if $L \in
\mathcal{K}_u^p(\mathbb{S}^n)$ such that $\gamma_u(K,L) \leq
\lambda$, then
\begin{equation} \label{anakin17}
L \subseteq \mathrm{conv}(K_{\lambda}) \qquad \mbox{and} \qquad K
\subseteq \mathrm{conv}(L_{\lambda}).
\end{equation}
We want to choose $\lambda = \overline{r}(K,\varepsilon)$ in such
a way that
\begin{equation} \label{padme17}
\mathrm{conv}(K_{\overline{r}(K,\varepsilon)}) \subseteq
K_{\varepsilon} \qquad \mbox{and} \qquad
\mathrm{conv}(L_{\overline{r}(K,\varepsilon)}) \subseteq
L_{\varepsilon}
\end{equation}
In order to compute $\overline{r}(K,\varepsilon)$ let $v, w \in
\mathbb{S}_u^+$ and denote by $J_v^w \in
\mathcal{K}_u^p(\mathbb{S}^n)$ the spherical segment connecting
$v$ and $w$. An elementary calculation shows that
\begin{equation} \label{chewy17}
\mathrm{conv}((J_v^w)_{\overline{c}(J_v^w,\varepsilon)})
\subseteq (J_v^w)_{\varepsilon},
\end{equation}
where
\[\overline{c}(J_v^w,\varepsilon)= \arcsin \left (\sin \varepsilon \cos\left (\frac{d(v,w)}{2} \right )    \right ).  \]
We define
\[\overline{r}(K,\varepsilon) = \min_{v,w \in K} \overline{c}(J_v^w,\varepsilon/2).  \]
Then, by (\ref{chewy17}),
\begin{eqnarray}
\mathrm{conv}(K_{\overline{r}(K,\varepsilon)})\!\!\!  & = & \!\!\!
\bigcup_{v,w \in K} \!
\mathrm{conv}((J_v^w)_{\overline{r}(K,\varepsilon)}) \subseteq
\!\! \bigcup_{v,w \in
K} \! \mathrm{conv}((J_v^w)_{\overline{c}(J_v^w,\varepsilon/2)}) \label{yoda42} \\
& \subseteq & \!\!\! \bigcup_{v,w \in K} \!
(J_v^w)_{\varepsilon/2} = K_{\varepsilon/2} \subseteq
K_{\varepsilon}. \label{yoda43}
\end{eqnarray}
This proves the first inclusion of (\ref{padme17}). To see the
second inclusion, note that
\[\overline{r}((J_v^w)_{\varepsilon/2},2\varepsilon) = \min_{v',w' \in (J_v^w)_{\varepsilon/2}} \overline{c}(J_{v'}^{w'},\varepsilon) \geq \overline{c}(J_v^w,\varepsilon/2).  \]
which follows from an elementary calculation and the fact, that
$\overline{c}(J_{v'}^{w'},\varepsilon)$ attains its minimum in
$(J_v^w)_{\varepsilon/2}$ when $d(v',w') = d(v,w)+\varepsilon$.
Thus, for $L \in \mathcal{K}_u^p(\mathbb{S}^n)$ such that
$\gamma_u(K,L) \leq \overline{r}(K,\varepsilon)$, it follows from
(\ref{anakin17}), (\ref{yoda42}), (\ref{yoda43}) that $L
\subseteq \mathrm{conv}(K_{\overline{r}(K,\varepsilon)})
\subseteq K_{\varepsilon/2}$ and we conclude
\[\overline{r}(L,2\varepsilon) \geq \!\! \min_{v,w \in K_{\varepsilon/2}} \!\! \overline{c}(J_v^w,\varepsilon) = \! \min_{v,w \in K} \overline{r}((J_v^w)_{\varepsilon/2},2\varepsilon)
\geq \! \min_{v,w \in K} \overline{c}(J_v^w,\varepsilon/2) =
\overline{r}(K,\varepsilon).\] Hence, using again (\ref{yoda42})
and (\ref{yoda43}), where $K$ is replaced by $L$,
\[\mathrm{conv}(L_{\overline{r}(K,\varepsilon)}) \subseteq \mathrm{conv}(L_{\overline{r}(L,2\varepsilon)}) \subseteq L_{\varepsilon}. \]
This proves the second inclusion of (\ref{padme17}) and, thus,
(\ref{ball1}). \hfill $\blacksquare$

\vspace{0.3cm}

Note that if $K, L \in \mathcal{K}_u^p(\mathbb{S}^n)$, then, by
(\ref{hausdsupp}) and Lemma 4.3,
\[\delta(g_u(K),g_u(L))= \max_{v\in \mathbb{S}_u} | \tan h_u(K,v)  - \tan h_u(L,v)|.\]
Thus, from Proposition \ref{schweregeb} and the continuity of the
tangent we obtain the following.

\begin{koro}\label{prop:gOhom}
The gnomonic projection is a homeomorphism between
$(\mathcal{K}_u^p(\mathbb{S}^n),\delta_s)$ and
$(\mathcal{K}(\mathbb{R}_u^n),\delta)$.
\end{koro}

Using Proposition \ref{prop:gOhom} and other basic properties of
the gnomonic projection, we can now prove the following
refinement of Theorem \ref{thm:2}.

\begin{theorem}\label{thm:2refined}
For every fixed $u \in \mathbb{S}^n$, the gnomonic projection $g_u$
induces a one-to-one correspondence between operations $*:
\mathcal{K}_u^p(\mathbb{S}^n) \times
\mathcal{K}_u^p(\mathbb{S}^n) \to \mathcal{K}_u^p(\mathbb{S}^n)$
which are $u$-projection covariant and operations $\overline{*}:
\mathcal{K}(\mathbb{R}_u^n) \times \mathcal{K}(\mathbb{R}_u^n)
\to \mathcal{K}(\mathbb{R}_u^n)$ \linebreak which are projection
covariant. Moreover, every such $u$-projection covariant operation
$*$ is continuous.

\end{theorem}
{\it Proof.} First assume that $*$ is $u$-projection covariant
and define an operation $\overline{*}:
\mathcal{K}(\mathbb{R}^n_u)\times\mathcal{K}(\mathbb{R}^n_u)\to\mathcal{K}(\mathbb{R}^n_u)$
by
\[K\, \overline{*}\, L = g_u( g_u^{-1}(K) * g_u^{-1}(L))\]
for $K, L \in \mathcal{K}(\mathbb{R}^n_u)$. Since for every
$k$-sphere $S$ containing $u$, there exists a linear subspace $V$
in $\mathbb{R}^n_u$ such that
\begin{equation*}
g_u^{-1}(K|V) = g_u^{-1}(K)|S
\end{equation*}
for all $K \in \mathcal{K}(\mathbb{R}_u^n)$, we obtain
\begin{align*}
(K|V)\,\overline{*} \,(L|V) &= g_u( g_u^{-1}(K|V)* g_u^{-1}(L|V)) = g_u((g_u^{-1}(K)|S) * (g_u^{-1}(L)|S))\\
        &= g_u((g_u^{-1}(K)*g_u^{-1}(L))|S) = g_u(g_u^{-1}(K)*g_u^{-1}(L))|V\\
        &= (K\,\overline{*}\, L)|V.
    \end{align*}
for all $K, L \in \mathcal{K}(\mathbb{R}_u^n)$. Thus,
$\overline{*}$ is projection covariant.

Now, let $\overline{*}:
\mathcal{K}(\mathbb{R}^n_u)\times\mathcal{K}(\mathbb{R}^n_u)\to\mathcal{K}(\mathbb{R}^n_u)$
be projection covariant and define $*:
\mathcal{K}_u^p(\mathbb{S}^n) \times
\mathcal{K}_u^p(\mathbb{S}^n) \to \mathcal{K}_u^p(\mathbb{S}^n)$
by
\[K * L = g_u^{-1}( g_u(K)\, \overline{*}\, g_u(L))\]
for $K, L \in \mathcal{K}_u^p(\mathbb{S}^n)$. Using a similar
argument as before, it is easy to show that $*$ is $u$-projection
covariant.

The continuity of an operation $*: \mathcal{K}_u^p(\mathbb{S}^n)
\times \mathcal{K}_u^p(\mathbb{S}^n) \to
\mathcal{K}_u^p(\mathbb{S}^n)$ which is $u$-projection covariant
is now a direct consequence of Theorem \ref{thm:gardner} and
Proposition \ref{prop:gOhom}. \hfill $\blacksquare$

\vspace{0.3cm}

Recall that the set $\mathcal{C} \subset
\mathcal{K}^p(\mathbb{S}^n) \times \mathcal{K}^p(\mathbb{S}^n)$
was defined by
\[ \mathcal{C}=\bigcup_{\substack{\O \in \S^n}}
\left( \mathcal{K}^p_u(\mathbb{S}^n) \times
\mathcal{K}^p_u(\mathbb{S}^n)  \right).\] By Theorem
\ref{thm:2refined}, the restriction of an operation $*:
\mathcal{C} \to \mathcal{K}^p(\mathbb{S}^n)$ which is projection
covariant to convex bodies contained in a fixed open hemisphere
is continuous. Therefore, we obtain:

\begin{koro} \label{cor:proj->cont} Every projection covariant operation $*\colon
\mathcal{C}\to\mathcal{K}^p(\mathbb{S}^n)$ is continuous.
\end{koro}

\vspace{1cm}

\centerline{\large{\bf{ \setcounter{abschnitt}{5}
\arabic{abschnitt}. Auxiliary results}}}

\reseteqn \alpheqn \setcounter{theorem}{0}

\vspace{0.6cm}

We continue in this section with our preparations for the proofs
of Theorems \ref{thm:1} and \ref{thm:3}. We prove three auxiliary
results which will be used at different stages in Section 6. We
begin by establishing first constraints on projection covariant
operations $*$ on $\mathcal{C}$.

\begin{lem}\label{lem:convexhull}
    If $*\colon\mathcal{C}\to\mathcal{K}^p(\mathbb{S}^n)$ is projection covariant, then either
    \begin{align}\label{cond:convexhull}
        K*L\subseteq \mathrm{conv}(K \cup L)
    \end{align}
    for all $(K,L)\in\mathcal{C}$   or
    \begin{align}
        K*L\subset -\mathrm{conv}(K \cup L)
    \end{align}
    for all $(K,L)\in\mathcal{C}$.
\end{lem}
{\it Proof.} For $u \in \mathbb{S}^n$, let $S_u$ denote the
$0$-sphere $\{-u, u\}$. By the projection covariance of $*$, we
have
\[(\{u\}*\{u\})|S_u = (\{u\}|S_u) * (\{u\}|S_u) = \{u\}*\{u\}.\]
Thus, $\{u\}*\{u\}\subseteq \{-u, u\}$. However, since
$\{u\}*\{u\} \in \mathcal{K}^p(\mathbb{S}^n)$, we must have either
$\{u\}*\{u\} = \{u\}$ or $\{u\}*\{u\}=\{-u\}$. Let
\[P = \{u \in \mathbb{S}^n: \{u\}*\{u\} = \{u\}\} \quad \mbox{and} \quad  N = \{u \in \mathbb{S}^n: \{u\}*\{u\} =
\{-u\}\}.\] Clearly, $P \cap N = \emptyset$ and $P \cup N=\S^n$.

Since, by Corollary \ref{cor:proj->cont}, $*$ is continuous, we
obtain for every sequence $u_i \in P$ with limit $u \in
\mathbb{S}^n$,
\[\{u\} * \{u\} =  \left\{ \lim u_i \right \} * \left\{ \lim u_i\right\} = \lim \left(\{u_i\}*\{u_i\}\right)
        = \lim \{u_i\} = \{u\}.\]
Thus, $u \in P$ which shows that $P$ is closed. In the same way,
we see that $N$ is closed. Consequently, we have either
$P=\mathbb{S}^n$ or $N=\S^n$.

First assume that $P = \mathbb{S}^n$ and let
$(K,L)\in\mathcal{C}$. Then there exists $u \in \mathbb{S}^n$ such
that $K$, $L \subset \mathbb{S}_u^+$ or, equivalently,
$\mathrm{conv}(K \cup L) \subset \mathbb{S}_u^+$. By the
projection covariance of $*$, we have
\[(K*L)|S_u = (K|S_u)*(L|S_u) = \{u\}*\{u\} = \{u\}.\]
Thus, $K*L\subset \mathbb{S}_u^+$ and we conclude that
\[K*L\subseteq \bigcap\{\mathbb{S}_u^+:u \in \mathbb{S}^n \text{ such that } \mathrm{conv}(K \cup L) \subset \mathbb{S}_u^+\}
        = \mathrm{conv}(K \cup L)\]
for all $(K,L)\in\mathcal{C}$.

Conversely, if $N = \mathbb{S}^n$, then we obtain $(K*L)|S_u =
\{-u\}$ and, therefore, $K*L \subset \mathbb{S}_u^- := -
\mathbb{S}_u^+$ whenever $\mathrm{conv}(K \cup L) \subset
\mathbb{S}_u^+$. This yields
    \begin{align*}
        K*L\subseteq \bigcap \{ \mathbb{S}_u^- : u\in\S^n \text{ such that } \mathrm{conv}(K \cup L) \subset \mathbb{S}_u^+\}
        = -\mathrm{conv}(K \cup L)
    \end{align*}
    for all $(K,L)\in\mathcal{C}$. \hfill $\blacksquare$

\vspace{0.3cm}

Our next lemma concerns spherical support functions of a spherical
segment contained in an open hemisphere.

\begin{lem} \label{eqn:Madd2}
For $u \in \mathbb{S}^n$, $v \in \mathbb{S}_u^+$, $w\in \mathbb{S}_u\cap \mathbb{S}_v$, and $-\frac{\pi}{2}< \alpha \leq \beta <
\frac{\pi}{2}$ let
\[I_u^w(\alpha,\beta) = \{u \cos \lambda  +w \sin \lambda : \lambda\in
[\alpha,\beta]\}.\] Then,
\[\tan h_v(I_u^w(\alpha,\beta),w) = \frac{\tan \beta}{u\cdot v} \quad \text{ and }
        \quad \tan h_v(I_u^w(\alpha,\beta),-w) = -\frac{\tan \alpha }{u\cdot v}.\]
\end{lem}
{\it Proof.} First note that by our definition of the spherical
support function
\[h_u(I_u^w(\alpha,\beta),w) = \beta \quad \mbox{and} \quad  h_u(I_u^w(\alpha,\beta),-w) = -\alpha.   \]
Let
\begin{align*}
A = g_v(I_u^w(\alpha,\alpha)) & = \frac{u \cos \alpha + w \sin
\alpha}{ (u\cdot v) \cos \alpha}-v,\\
B=g_v(I_u^w(\beta,\beta)) &= \frac{u \cos \beta + w \sin
\beta}{(u\cdot v)\cos \beta }-v.
\end{align*}
By Lemma \ref{gnombasic1} (b), $g_v(I_u^w(\alpha,\beta))$ is the
line segment in $\mathbb{R}^n_v$ in direction $w$ with endpoints
$A$ and $B$. Thus, by Lemma \ref{suppsupp} and the definition of
(Euclidean) support functions, we obtain
\begin{align*}
\tan h_v(I_u^w(\alpha,\beta),w) &= h(g_v(I_u^w(\alpha,\beta)),w)
= w \cdot B = \frac{\tan\beta}{u\cdot v},\\
\tan h_v(I_u^w(\alpha,\beta),-w) &=
h(g_v(I_u^w(\alpha,\beta)),-w)= -w \cdot A
        = -\frac{\tan \alpha}{u\cdot v}.
    \end{align*}

\vspace{-0.4cm}

\hfill $\blacksquare$

\vspace{0.4cm}

In view of Lemma \ref{lem:convexhull}, Theorem \ref{thm:2refined},
and Theorem \ref{thm:gardner}, the following result will be
useful in the proof of Theorem \ref{thm:3}.

\begin{lem}\label{lem:Mconvexhull}
Let $M \subseteq \mathbb{R}^4$ be closed and convex. If for all
$a, b, c, d \in \mathbb{R}$ such that $-a \leq b$ and $-c \leq d$,
\begin{align}\label{cond:mplus1}
 h_M(a,b,c,d) \leq \max\{b, d\},
\end{align}
then
\begin{align*}
M\subseteq
\{(\lambda_2,\lambda_1+\lambda_2,\lambda_3,1-\lambda_1+\lambda_3)
\in \mathbb{R}^4 : \lambda_1\in[0,1],\lambda_2\leq 0,\lambda_3\leq
0\}.
\end{align*}
\end{lem}
{\it Proof.} For $z=(-1,1,-1,1)$, we obtain from
(\ref{cond:mplus1}) that
\[h(M,z) \leq 1 \qquad \mbox{and} \qquad h(M,-z) \leq -1.\]
Since $-h(M,-z) \leq h(M,z)$, we conclude that $-h(M,-z)=h(M,z) =
1$ or, equivalently,
\begin{align}\label{step1}
   M\subseteq \{x\in\R^4 : -x_1+x_2-x_3+x_4 = 1\}.
\end{align}

By (\ref{cond:mplus1}), we also have $h_M(1,0,0,0)\leq 0$
    and $h_M(0,0,1,0)\leq 0$. Thus,
    \begin{align}\label{step2}
        M \subseteq \{x\in\R^4: x_1 \leq 0, x_3\leq 0\}.
    \end{align}

Finally, we deduce from (\ref{cond:mplus1}) that
    \begin{align*}
        h_M(-1,1,0,0)\leq 1 \qquad \text{ and } \qquad h_M(1,-1,0,0)\leq
        0,
    \end{align*}
    as well as
    \begin{align*}
        h_M(0,0,-1,1) \leq 1 \qquad \text{ and } \qquad h_M(0,0,1,-1) \leq 0.
    \end{align*}
Consequently,
    \begin{align}\label{step3}
        M\subseteq\{x\in\R^4: 0\leq x_2-x_1\leq 1 \text{ and } 0\leq x_4-x_3 \leq 1\}.
    \end{align}

    Combining (\ref{step1}), (\ref{step2}), and (\ref{step3}),
    completes the proof.
\hfill $\blacksquare$

\vspace{0.4cm}

The importance for us of the set
\[E:= \{(\lambda_2,\lambda_1+\lambda_2,\lambda_3,1-\lambda_1+\lambda_3) \in \mathbb{R}^4:
\lambda_1\in[0,1],\lambda_2\leq 0,\lambda_3\leq 0\} \] follows
from
\[h_E(h_{-K}(x),h_K(x),h_{-L}(x),h_L(x)) = h_{\mathrm{conv}(K \cup L)}(x).\]

\vspace{1cm}

\centerline{\large{\bf{ \setcounter{abschnitt}{6}
\arabic{abschnitt}. Proofs of Theorems \ref{thm:1} and
\ref{thm:3}}}}

\reseteqn \alpheqn \setcounter{theorem}{0}

\vspace{0.6cm}

After these preparations, we are now in a position to first proof
Theorem~\ref{thm:3} and then complete the proof of Theorem
\ref{thm:1}. In order to enhance the readability of several
formulas below, we write $\tan(x_1,\ldots,x_k)$ for the vector
$(\tan x_1,\ldots,\tan x_k)$ and $\arctan(x_1,\ldots,x_k)$ is
defined similarly.

\begin{theorem} \label{langabercool}
An operation $*\colon \mathcal{C} \to
\mathcal{K}^p(\mathbb{S}^n)$ is projection covariant if and only
if it is either  $K * L = \mathrm{conv}(K\cup L)$ or $K * L = -
\mathrm{conv}(K \cup L)$ for all $(K,L) \in \mathcal{C}$ or it is
trivial, that is, $K * L = K$, or $K * L = -K$, or $K * L = L$, or
$K * L = -L$ for all $(K,L) \in \mathcal{C}$.
\end{theorem}
{\it Proof.} By Lemma \ref{lem:convexhull}, we may assume that
\begin{equation} \label{enthalten17}
K * L \subseteq \mathrm{conv}(K \cup L)
\end{equation}
holds for all $(K,L) \in \mathcal{C}$ (otherwise, replace $*$ by
$*^-: \mathcal{C} \to \mathcal{K}^p(\mathbb{S}^n)$ defined by
$K*^-L = -(K*L)$). In particular, for every $u \in \mathbb{S}^n$,
the range of the restriction of $*$ to
$\mathcal{K}^p_u(\mathbb{S}^n) \times
\mathcal{K}^p_u(\mathbb{S}^n)$ lies in
$\mathcal{K}^p_u(\mathbb{S}^n)$.

In the proof of Theorem \ref{thm:2refined} we have seen that, for
every $u \in \mathbb{S}^n$, there exists a (unique) projection
covariant operation $\overline{*}_u\colon
\mathcal{K}(\mathbb{R}^n_u)\times \mathcal{K}(\mathbb{R}^n_u) \to
\mathcal{K}(\mathbb{R}^n_u)$ such that
\begin{align*}
\overline{K}\,\,\overline{*}_u\, \overline{L} = g_u(
g_u^{-1}(\overline{K}) * g_u^{-1}(\overline{L}))
\end{align*}
for all $\overline{K}$, $\overline{L}
\in\mathcal{K}(\mathbb{R}^n_u)$. Thus, by Theorem
\ref{thm:gardner}, there exists a nonempty closed convex set $M_u
\subset \R^4$ such that
\[h_{\overline{K}\,\,\overline{*}_u \overline{L}}(v) =
h_{M_u}(h_{\overline{K}}(-v),h_{\overline{K}}(v),h_{\overline{L}}(-v),h_{\overline{L}}(v))\]
for all $v \in \mathbb{S}_u$. Therefore, Lemma \ref{suppsupp}
yields
\begin{align}
\tan h_u(K*L,v) & = h_{g_u(K*L)}(v)  = h_{g_u(K)\,\overline{*}_u\, g_u(L)}(v) \label{bw117}\\
 & = h_{M_u}\!(h_{g_u(K)}(-v),h_{g_u(K)}(v),h_{g_u(L)}(-v),h_{g_u(L)}(v)) \nonumber \\
 & =
 h_{M_u}\!(\tan(h_u(K,-v),h_u(K,v),h_u(L,-v),h_u(L,v)))
 \label{bw442}
\end{align}
for all $K,L\in\mathcal{K}_u^p(\mathbb{S}^n)$. Thus, since
$-h_{\overline{K}}(-v) \leq h_{\overline{K}}(v)$ for every
$\overline{K} \in \mathcal{K}(\mathbb{R}^n_u)$ and every $v \in
\mathbb{S}_u$, the restriction of $*$ to
$\mathcal{K}_u^p(\mathbb{S}^n) \times
\mathcal{K}_u^p(\mathbb{S}^n)$ is completely determined by the
values $h_{M_u}(a,b,c,d)$, where $-a \leq b$ and $-c \leq d$.
Next, we want to show that for such $a,b,c,d \in \mathbb{R}$,
\begin{equation} \label{abcd}
h_{M_u}(a,b,c,d) = h_{M_v}(a,b,c,d)
\end{equation}
whenever $v \in \mathbb{S}_u^+$. To this end, let
$-\frac{\pi}{2}< \alpha \leq \beta < \frac{\pi}{2}$ and
$-\frac{\pi}{2}< \varphi \leq \psi < \frac{\pi}{2}$. For every $u
\in \mathbb{S}^n$ and $w \in \mathbb{S}_u$, the $u$-projection
covariance of $*$ implies that there exist $\sigma, \tau$ such
that $-\frac{\pi}{2}< \sigma \leq \tau < \frac{\pi}{2}$ and
\begin{align} \label{bwyoda1}
        I_u^w(\alpha,\beta)*I_u^w(\varphi,\psi) =
        I_u^w(\sigma,\tau),
\end{align}
where we have used the notation from Lemma \ref{eqn:Madd2}
 for spherical segments $I_u^w$. Since, for $-\frac{\pi}{2}< \xi \leq \zeta <
 \frac{\pi}{2}$, we have $h_u(I_u^w(\xi,\zeta),-w) = -\xi$ and
 $h_u(I_u^w(\xi,\zeta),w) = \zeta$, we obtain on the one hand from
 (\ref{bwyoda1}), (\ref{bw117}), and (\ref{bw442}),
\begin{align*}
\tan \tau  & = \tan h_u(I_u^w(\sigma,\tau),w) = \tan h_u(I_u^w(\alpha,\beta)*I_u^w(\varphi,\psi),w)\\
           &=h_{M_u}(\tan(-\alpha,\beta,-\varphi,\psi)).
    \end{align*}
For $v \in \mathbb{S}_u^+$ and $w \in \mathbb{S}_u \cap
\mathbb{S}_v$, we obtain from Lemma \ref{eqn:Madd2} and again
(\ref{bwyoda1}), (\ref{bw117}), and (\ref{bw442}),
\begin{align*}
\tan \tau & = (u \cdot v) \tan h_v(I_u^w(\sigma,\tau),w) = (u\cdot v) \tan h_v(I_u^w(\alpha,\beta)*I_u^w(\varphi,\psi),w)\\
          & = (u \cdot v) h_{M_v}\left(\frac{\tan(-\alpha,\beta,-\varphi,\psi)}{u\cdot v}
          \right) = h_{M_v}(\tan(-\alpha,\beta,-\varphi,\psi))
\end{align*}
which proves (\ref{abcd}). Since $u \in \mathbb{S}^n$, $v \in
\mathbb{S}_u^+$, and $\alpha, \beta, \varphi, \psi$ were
arbitrary, we conclude from (\ref{bw117}), (\ref{bw442}), and
(\ref{abcd}) that there exists a nonempty closed convex set $M
\subseteq \mathbb{R}^4$, \emph{independent} of $u \in
\mathbb{S}^n$, such that
\begin{equation} \label{superformel}
\tan h_u(K*L,v) =
h_{M}(\tan(h_u(K,-v),h_u(K,v),h_u(L,-v),h_u(L,v)))
\end{equation}
for all $K, L \in \mathcal{K}_u^p(\mathbb{S}^n)$ and $v \in
\mathbb{S}_u$.

\vspace{0.2cm}

To complete the proof, we have to show that for $-a \leq b$ and
$-c \leq d$, the support function $h_M$ satisfies one of the
following three conditions:
\begin{enumerate}
\item[(i)] $h_M(a,b,c,d) = b$, that is, $K * L = K$ for $(K,L)
\in \mathcal{C}$;
\item[(ii)] $h_M(a,b,c,d) = d$, that is, $K * L = L$ for $(K,L)
\in \mathcal{C}$;
\item[(iii)] $h_M(a,b,c,d) =\max\{b,d\}$, that is, $K * L = \mathrm{conv}(K \cup L)$ for $(K,L)
\in \mathcal{C}$.
\end{enumerate}

From (\ref{enthalten17}) and (\ref{superformel}), we deduce that
\begin{align} \label{maxent}
h_M(a,b,c,d)\leq \max\{b,d\}
\end{align}
whenever $-a \leq b$ and $-c \leq d$. Moreover, since $-h_u(K *
L,-v) \leq h_u(K * L,v)$ for all $K, L \in
\mathcal{K}_u^p(\mathbb{S}^n)$ and $v \in \mathbb{S}_u$, we
deduce from (\ref{superformel}) that
\begin{align}\label{eqn:hMineq}
-h_M(b,a,d,c) \leq h_M(a,b,c,d).
\end{align}

Next, we want to show that for all $-\frac{\pi}{2}<\alpha\leq
\beta < \frac{\pi}{2}$, $-\frac{\pi}{2}< \varphi \leq \psi <
\frac{\pi}{2}$, and $-\frac{\pi}{2}+\max\{\beta,\psi\}< \eta <
\frac{\pi}{2}+\min\{\alpha,\varphi\}$, we have
\begin{align}\label{prop:rotation}
\arctan h_M(\tan \Lambda) = \arctan h_M(\tan(\Lambda + \Theta))+
\eta,
    \end{align}
where $\Lambda = (-\alpha,\beta,-\varphi,\psi)$ and $\Theta
=(\eta,-\eta,\eta,-\eta)$. In order to prove
(\ref{prop:rotation}), let $u \in \S^n$, $v \in \mathbb{S}_u$ and
define
\begin{align*}
u' = u \cos \eta  - v \sin \eta \qquad \mbox{and} \qquad v' = v
\cos \eta + u \sin \eta.
    \end{align*}
Note that $u'$ and $v'$ are rotations of $u$ and $v$ in the plane
$\mathrm{span}\{u, v\}$ by an angle $-\eta$. Therefore, for every
$\lambda \in [0,2\pi)$,
\begin{align*}
u' \cos \lambda + v' \sin \lambda  = u \cos(\lambda-\eta) + v
\sin(\lambda-\eta).
\end{align*}
Hence,
\begin{align} \label{luke42}
I_ {u'}^{v'}(\alpha,\beta) = I_u^v(\alpha-\eta,\beta-\eta)
\subseteq \mathbb{S}_u^+.
\end{align}
Now, let
    \begin{align*}
        \sigma  & = - h_u(I_u^v(\alpha-\eta,\beta-\eta)*I_u^v(\varphi-\eta,\psi-\eta),-v), \\
        \tau & =
        h_u(I_u^v(\alpha-\eta,\beta-\eta)*I_u^v(\varphi-\eta,\psi-\eta),v),
    \end{align*}
and
\begin{align*}
        \sigma' & = -h_{u'}(I_{u'}^{v'}(\alpha,\beta) * I_{u'}^{v'}(\varphi,\psi),-v'),\\
        \tau'   & = h_{u'}(I_{u'}^{v'}(\alpha,\beta) * I_{u'}^{v'}(\varphi,\psi),v').
\end{align*}
By the $u$-projection covariance and the $u'$-projection
covariance of $*$ and (\ref{luke42}), we obtain
\begin{align*}
I_{u'}^{v'}(\sigma',\tau')&=I_{u'}^{v'}(\alpha,\beta) * I_{u'}^{v'}(\varphi,\psi) = I_u^v(\alpha-\eta,\beta-\eta)*I_u^v(\varphi-\eta,\psi-\eta)\\
        &= I_u^v(\sigma,\tau) = I_{u'}^{v'}(\sigma+\eta,\tau+\eta).
    \end{align*}
Thus, $\tau' = \tau + \eta$. Using (\ref{superformel}) and the
definitions of $\tau$ and $\tau'$, we obtain
(\ref{prop:rotation}).

\pagebreak

From applications of (\ref{prop:rotation}) with $\Lambda = \pm
(-\alpha,\alpha,\alpha,-\alpha)$ and $\eta = \pm \alpha$, where
$\alpha\in[0,\frac{\pi}{4})$, we obtain
\begin{eqnarray}
\arctan(h_M(-1,1,1,-1)\tan \alpha) = \arctan(h_M(0,0,1,-1)\tan(2\alpha))+\alpha, \phantom{W} & & \label{hansolo42}\\
\arctan(h_M(-1,1,1,-1)\tan \alpha) =
\arctan(h_M(-1,1,0,0)\tan(2\alpha))-\alpha, \phantom{W} & &
\label{leia42}
\end{eqnarray}
and
\begin{eqnarray}
\arctan(h_M(1,-1,-1,1)\tan \alpha) = \arctan(h_M(0,0,-1,1)\tan(2\alpha))-\alpha, \phantom{W} & & \label{obiwan42}\\
\arctan(h_M(1,-1,-1,1)\tan \alpha) =
\arctan(h_M(1,-1,0,0)\tan(2\alpha))+\alpha. \phantom{\,W} & &
\label{chewy42}
\end{eqnarray}

On the one hand, using (\ref{hansolo42}) and (\ref{leia42}), it is
not difficult to show that either
\begin{eqnarray} \label{cond:1}
h_M(-1,1,1,-1)=1, \quad h_M(0,0,1,-1)=0, \quad h_M(-1,1,0,0)=1,
\phantom{WW} & &
\end{eqnarray}
or
\begin{eqnarray} \label{cond:3}
h_M(-1,1,1,-1)=-1, \quad h_M(0,0,1,-1)=-1, \quad h_M(-1,1,0,0)=0.
\phantom{W} & &
\end{eqnarray}
On the other hand, by (\ref{obiwan42}) and (\ref{chewy42}), we
have either
\begin{eqnarray} \label{cond:2}
h_M(1,-1,-1,1)=1, \quad h_M(0,0,-1,1)=1, \quad h_M(1,-1,0,0)=0,
\phantom{W} & &
\end{eqnarray}
or
\begin{eqnarray} \label{cond:4}
h_M(1,-1,-1,1)=-1, \quad h_M(0,0,-1,1)=0, \quad h_M(1,-1,0,0)=-1.
\phantom{W} & &
\end{eqnarray}
Note that, since $-h_M(1,-1,-1,1) \leq h_M(-1,1,1,-1)$,
(\ref{cond:3}) and (\ref{cond:4}) cannot both be satisfied. Also
recall that by Lemma \ref{lem:Mconvexhull}, we have
    \begin{align*}
        M\subseteq E=\{(\lambda_2,\lambda_1+\lambda_2,\lambda_3,1-\lambda_1+\lambda_3):
            \lambda_1\in[0,1], \lambda_2,\lambda_3\leq 0\}.
    \end{align*}
and let
\begin{align*}
E_0&=\{(\lambda_2,\lambda_2,\lambda_3,1+\lambda_3) : \lambda_2
\leq 0,\lambda_3\leq 0\},\\
E_1&=\{(\lambda_2,1+\lambda_2,\lambda_3,\lambda_3) : \lambda_2
\leq 0, \lambda_3\leq 0\}.
    \end{align*}

If (\ref{cond:1}) holds, then $h_M(-1,1,0,0) = 1$ and, since $M
\subseteq E$, we have
\begin{align*}
1=\max\{\lambda_1 \in [0,1]
:(\lambda_2,\lambda_1+\lambda_2,\lambda_3,1-\lambda_1+\lambda_3)\in
M\}.
\end{align*}
Thus, there are $\lambda_2$, $\lambda_3 \leq 0$, such that
$(\lambda_2,1+\lambda_2,\lambda_3,\lambda_3)\in M$ or,
equivalently, $M \cap E_1$ is nonempty. Similarly, it follows from
(\ref{cond:2}) that $M\cap E_0$ is nonempty.

\pagebreak

If (\ref{cond:3}) holds, we have $h_M(-1,1,0,0) = 0$ and we deduce
that
\begin{align*}
0 = \max\{\lambda_1 \in [0,1]:
(\lambda_2,\lambda_1+\lambda_2,\lambda_3,1-\lambda_1+\lambda_3)\in
M\}
\end{align*}
which yields $M \subseteq E_0$. Analogously, (\ref{cond:4})
implies $M\subseteq E_1$.

Next, an application of (\ref{prop:rotation}) with $\Lambda =
(0,\alpha,0,\alpha)$ and $\eta = \alpha$, where again
$\alpha\in[0,\frac{\pi}{4})$, yields
\begin{align*}
\arctan(h_M(0,1,0,1)\tan \alpha) = \arctan(h_M(1,0,1,0) \tan
\alpha) + \alpha
\end{align*}
Clearly, this is possible if and only if either
\begin{align} \label{cond:5}
h_M(0,1,0,1) = 1, \quad  h_M(1,0,1,0)=0,
\end{align}
or
\begin{align} \label{cond:6}
h_M(0,1,0,1) = 0, \quad  h_M(1,0,1,0)=-1.
\end{align}
However, (\ref{cond:6}) contradicts (\ref{eqn:hMineq}) and is
therefore not possible.

From (\ref{cond:5}) and the fact that $M \subseteq E$, we infer
\begin{align*}
0 = \max \{\lambda_2+\lambda_3 : \lambda_2, \lambda_3 \leq 0
\mbox{ and }
(\lambda_2,\lambda_1+\lambda_2,\lambda_3,1-\lambda_1+\lambda_3)\in
M\}
\end{align*}
which implies
    \begin{align}\label{cond:bla}
        M\cap \{(0,\lambda_1,0,1-\lambda_1): \lambda_1\in[0,1]\} \neq
        \emptyset.
    \end{align}

\vspace{0.2cm}

For the final part of the proof, we distinguish three cases:
\begin{itemize}
\item[(i)] (\ref{cond:1}) and (\ref{cond:4}) hold, in particular, $M
\subseteq E_1$;
\item[(ii)] (\ref{cond:3}) and (\ref{cond:2}) hold, in particular, $M \subseteq
E_0$;
\item[(iii)] (\ref{cond:1}) and (\ref{cond:2}) hold.
\end{itemize}

In case (i), $M \subseteq E_1$ and (\ref{cond:bla}) imply that
$e_2 \in M$. Using (\ref{maxent}), we conclude that $h_M(a,b,c,d)
=b$, that is, $K
* L = K$ for $(K,L) \in \mathcal{C}$.

Similarly, in case (ii), $M \subseteq E_0$ and (\ref{cond:bla})
imply that $e_4 \in M$. Using again (\ref{maxent}), we obtain
$h_M(a,b,c,d) = d$, that is, $K * L = L$ for $(K,L) \in
\mathcal{C}$.

It remains to show that in case (iii), we have $e_2$, $e_4\in M$
which, by (\ref{maxent}), implies that $h_M(a,b,c,d) =
\max\{b,d\}$ or $K * L = \mathrm{conv}(K \cup L)$ for $(K,L) \in
\mathcal{C}$. \linebreak To this end, we apply again
(\ref{prop:rotation}) with $\Lambda = (0,\alpha,0,0)$ and $\eta =
\alpha$, where $\alpha \in [0,\frac{\pi}{4})$ to obtain
\begin{align*}
\arctan(h_M(0,1,0,0)\tan \alpha) &= \arctan(h_M(1,0,1,-1) \tan
\alpha) + \alpha
\end{align*}
This is possible if and only if either
\begin{align} \label{cond:7}
h_M(0,1,0,0) = 1, \quad  h_M(1,0,1,-1) = 0,
\end{align}
or
\begin{align} \label{cond:9}
h_M(0,1,0,0) = 0, \quad h_M(1,0,1,-1) = -1.
\end{align}

Assume that (\ref{cond:9}) holds. Then, by (\ref{maxent}),
(\ref{eqn:hMineq}), and the subadditivity of $h_M$, we obtain
\begin{align}\label{eqn:proofstep}
-1 \leq h_M(1,1,1,-1) \leq h_M(1,0,1,-1) + h_M(0,1,0,0) = -1
\end{align}
Hence, $h_M(1,1,1,-1) = -1$.

Now, consider the convex bodies $\overline{K}=[-e_2,e_2]$ and
$\overline{L}=\{e_1\}$ in $\mathbb{R}^n$. Then,
$h_{\overline{K}}(x) = |e_2 \cdot x|$ and $h_{\overline{L}}(x) =
e_1 \cdot x$ for $x \in \mathbb{R}^n$, and we obtain from
(\ref{cond:1}) and $h_M(1,1,1,-1) = -1$,
\begin{align*}
h_M(h_{\overline{K}}(e_1),h_{\overline{K}}(-e_1),h_{\overline{L}}(e_1),h_{\overline{L}}(-e_1))&=h_M(0,0,1,\!-1)
= 0,\\
h_M(h_{\overline{K}}(e_1\!\!+\!e_2),h_{\overline{K}}(-e_1\!\!-\!e_2),h_{\overline{L}}(e_1\!\!+\!e_2),h_{\overline{L}}(-e_1\!\!-\!e_2))&= h_M(1,1,1,\!-1) = -1,\\
h_M(h_{\overline{K}}(e_1\!\!-\!e_2),h_{\overline{K}}(e_2\!\!-\!e_1),h_{\overline{L}}(e_1\!\!-\!e_2),h_{\overline{L}}(e_2\!\!-\!e_1))
& = h_M(1,1,1,\!-1) = -1.
\end{align*}
Since
$h_M(h_{-\overline{K}},h_{\overline{K}},h_{-\overline{L}},h_{\overline{L}})$
defines a support function of a convex body $\overline{Z}$ in
$\mathbb{R}^n$, we infer
    \begin{align*}
        0 = h_{\overline{Z}}(-2e_1) \geq h_{\overline{Z}}(-e_1-e_2) + h_{\overline{Z}}(-e_1+e_2) = -2
    \end{align*}
which contradicts the subadditivity of $h_{\overline{Z}}$. Thus,
(\ref{cond:9}) cannot hold.

Another application of (\ref{prop:rotation}) with $\Lambda =
(0,\alpha,\alpha,0)$ and $\eta = -\alpha$, where $\alpha \in
[0,\frac{\pi}{4})$, yields
\begin{align*}
\arctan(h_M(0,1,1,0)\tan \alpha) &=
\arctan(h_M(\tan(-\alpha,2\alpha,0,\alpha))) - \alpha.
\end{align*}
Consequently,
\begin{align*}
h_M\left(-\frac{\tan \alpha}{\tan(2\alpha)},1,0,\frac{\tan
\alpha}{\tan(2\alpha)}\right)=\frac{\tan(\arctan(h_M(0,1,1,0)\tan
\alpha)+\alpha)}{\tan(2\alpha)}.
    \end{align*}
By letting $\alpha\to \frac{\pi}{4}$ and using (\ref{cond:7}), we
deduce that $h_M(0,1,1,0) = 1$. Since $M \subseteq E$, this yields
\begin{align*}
1=\max\{ \lambda_1+\lambda_2+\lambda_3 :
       (\lambda_2,\lambda_1+\lambda_2,\lambda_3,(1-\lambda_1)+\lambda_3)\in M\}
    \end{align*}
which, in turn, implies that $e_2 \in M$.

The proof that $e_4 \in M$ is now very similar. We first use
(\ref{prop:rotation}) with $\Lambda = (0,0,0,\alpha)$ and $\eta =
\alpha$ to deduce that
\begin{align*}
h_M(0,0,0,1) = 1, \quad h_M(1,-1,1,0) = 0.
\end{align*}
Using this and another application of (\ref{prop:rotation}) with
$\Lambda = (\alpha,0,0,\alpha)$ and $\eta = -\alpha$, finally
leads to $h_M(1,0,0,1) =1$. From this and $M \subseteq E$,
follows $e_4 \in M$ which completes the proof.
\hfill$\blacksquare$

\vspace{0.3cm}

Using Theorem \ref{langabercool}, we can now also complete the
proof of Theorem \ref{thm:1}:

\begin{theorem} An operation $*\colon \mathcal{K}^p(\mathbb{S}^n) \times \mathcal{K}^p(\mathbb{S}^n) \to \mathcal{K}^p(\mathbb{S}^n)$
is projection covariant and continuous if and only if either $K *
L = K$, or $K * L = -K$, or $K * L = L$, or $K * L = -L$ for all
$K, L \in \mathcal{K}^p(\mathbb{S}^n)$.
\end{theorem}
{\it Proof.} By Theorem \ref{langabercool}, it is sufficient to
prove that the convex hull does not admit a continuous extension
to a map from $\K^p(\S^n)\times\K^p(\S^n)$ to $\K^p(\S^n)$. In
order to show this, let $u\in \mathbb{S}^n$, $v\in\mathbb{S}_u$,
and consider the spherical segments $K = I_u^v(-\frac{\pi}{2},0)$
and $L_\varepsilon = I_u^v(0,\frac{\pi}{2}-\varepsilon)$, where
$\varepsilon > 0$. Then $(K,L_\varepsilon)\in\mathcal{C}$
converges in the Hausdorff metric to
$(K,L_0)\in\K^p(\S^n)\times\K^p(\S^n)$ as $\varepsilon
\rightarrow 0^+$. However,
\begin{align*}
\lim_{\varepsilon\to 0^+} \mathrm{conv}(K \cup L_{\varepsilon}) =
\lim_{\varepsilon\to 0^+}
I_u^v\left(-\frac{\pi}{2},\frac{\pi}{2}-\varepsilon\right)
        =I_u^v\left(-\frac{\pi}{2},\frac{\pi}{2}\right) \not \in \K^p(\S^n).
    \end{align*}

\vspace{-0.5cm}

\hfill$\blacksquare$

\vspace{0.4cm}

We remark, that it is also not difficult to show that the convex hull is {\it not} continuous as a map
from $\K^p(\S^n)\times\K^p(\S^n)$ to $\mathcal{K}(\mathbb{S}^n)$.

\vspace{1cm}

\centerline{\large{\bf{ \setcounter{abschnitt}{7}
\arabic{abschnitt}. Section covariant operations}}}

\reseteqn \alpheqn \setcounter{theorem}{0}

\vspace{0.6cm}

In this final section, first we briefly recall a characterization of rotation and section covariant operations between Euclidean star sets established in \textbf{\cite{GardHugWeil2013}}.
Than, we discuss basic properties of spherical star sets in order to eventually prove a corresponding result to Theorem \ref{thm:2} for rotation and section covariant operations between them.

\pagebreak

A subset $L$ of $\mathbb{R}^n$ is called \emph{star-shaped} with respect to $o$ if every line through the origin intersects $L$ in a (possibly degenerate) closed line segment.
A \emph{star set} in $\mathbb{R}^n$ is a compact set that is star-shaped with respect to $o$.
The radial function $\rho(L,\cdot) = \rho_L: \mathbb{R}^n \backslash \{o\} \rightarrow [0,\infty)$ of a star set $L$ is defined by
\[\rho(L,x) = \max\{\lambda \geq 0: \lambda x \in L\}, \qquad x \in \mathbb{R}^n \backslash \{o\}.  \]
Radial functions are $-1$-homogeneous, that is, $\rho(L,\lambda x) = \lambda^{-1}\rho(L,x)$ for all
$x \in \mathbb{R}^n \backslash \{o\}$ and $\lambda > 0$, and are therefore often regarded as functions on $\mathbb{S}^{n-1}$. If $\rho(L,\cdot)$ is
positive and continuous, we call $L$ a \emph{star body}. If $K \in \mathcal{K}(\mathbb{R}^n)$ contains the origin in its interior,
then $K$ is a star body and we have
\begin{equation} \label{polar17}
\rho(K^*,\cdot) = \frac{1}{h(K,\cdot)} \qquad \mbox{and} \qquad h(K^*,\cdot) = \frac{1}{\rho(K,\cdot)},
\end{equation}
where $K^*$ denotes the \emph{polar body} of $K$ defined by
\[K^* = \{x \in \mathbb{R}^n: x \cdot y \leq 1 \mbox{ for all } y \in K\}.  \]

The \emph{radial distance} $\widetilde{\delta}(K,L)$ between two star sets $K$ and $L$ in $\mathbb{R}^n$ is defined by
\begin{equation} \label{defraddist}
\widetilde{\delta}(K,L) = \|\rho(K,\cdot) - \rho(L,\cdot)\|_{\infty}.
\end{equation}

We denote by $\mathcal{S}(\mathbb{R}^n)$ the space of all star sets in $\mathbb{R}^n$ endowed with the radial distance.
The {\it radial sum} $K\,\widetilde{+}\,L$ of $K$, $L \in \mathcal{S}(\mathbb{R}^n)$ is defined as the star set such that
\[\rho(K\,\widetilde{+}\,L,\cdot) = \rho(K,\cdot) + \rho(L,\cdot).  \]
More generally, for any $p > 0$, the \emph{$L_p$ radial sum} $K\,\widetilde{+}_p\,L$ of $K$, $L \in \mathcal{S}(\mathbb{R}^n)$ is defined by
\[\rho(K\,\widetilde{+}_p\,L,\cdot)^p = \rho(K,\cdot)^p + \rho(L,\cdot)^p.   \]
Lutwak \textbf{\cite{Lutwak:1996}} showed that in the same way as the $L_p$ Minkowski addition leads to the $L_p$ Brunn--Minkowski theory, $L_p$ radial addition leads to a dual $L_p$ Brunn--Minkowski theory (see also \textbf{\cite{Gardner:2006}} and the references therein).

While $L_p$ radial addition is \emph{not} projection covariant, the $L_p$ radial sum of star sets is \emph{section covariant}, that is,
\[(K \cap V) \,\widetilde{+}_p\, (L \cap V) = (K\,\widetilde{+}_p\,L) \cap V  \]
for every linear subspace $V$ of $\mathbb{R}^n$. It is also $\mathrm{GL}(n)$ covariant and therefore, in particular,
covariant with respect to rotations.

A complete classification of all rotation and section covariant binary operations between
star sets in $\mathbb{R}^n$ was established by Gardner, Hug, and Weil and can be stated as follows:

\pagebreak

\begin{theorem} \label{thm:starshaped} {\bf (\hspace{-0.15cm} \cite{GardHugWeil2013})}
An operation $*: \mathcal{S}(\R^n)\times \mathcal{S}(\R^n) \to \mathcal{S}(\R^n)$ is rotation and section covariant
if and only if there exists a function $f:[0,\infty)^4\to \mathbb{R}$ such that, for all $K$, $L\in\mathcal{S}(\R^n)$,
\begin{align*}
\rho_{K*L}(v) = f(\rho_{-K}(v),\rho_K(v),\rho_{-L}(v),\rho_L(v)), \qquad v \in \mathbb{S}^{n-1}.
\end{align*}
\end{theorem}

\vspace{0.2cm}

We turn now to star sets in $\mathbb{S}^n$. We call a subset $L$ of $\mathbb{S}^n$ a {\it (spherical) star set}
with respect to $u \in L$ if $L \cap S_{u,v}$ is a (possibly degenerate) closed spherical segment for all $v\in \S_u$.
We denote by $\mathcal{S}_u(\S^n)$ the class of all spherical star sets with respect to $u$ and we write
$\mathcal{S}^p_u(\S^n)$ for the subclass of \emph{proper star sets} with respect
to $u$, that is, star sets with respect to $u$ contained in $\S_u^+$.

\vspace{0.3cm}

\noindent {\bf Definition} For $u \in \mathbb{S}^n$ and a proper star set $L \in
\mathcal{S}_u^p(\mathbb{S}^n)$, the \emph{spherical radial
function} $\rho_u(L,\cdot): \S_u \to [0,\frac{\pi}{2})$ of $L$ is defined by
\[\rho_u(L,v) = \max\{\alpha \geq 0: u \cos \alpha  + v \sin \alpha \in L\}.\]

\vspace{0.2cm}

Note that, for every $v \in \mathbb{S}_u$, we have
\begin{align*}
u\cos \rho_u(L,v) + v\sin \rho_u(L,v) \in \partial L.
\end{align*}

The counterparts to Lemma \ref{gnombasic1} (c) and Lemma \ref{suppsupp} in the setting of spherical star sets are the contents of our next lemma.

\begin{lem} \label{radrad}
For $u \in \mathbb{S}^n$, the following statements hold:
\begin{enumerate}
\item[(a)] The gnomonic projection maps $\mathcal{S}_u^p(\mathbb{S}^n)$ bijectively to $\mathcal{S}(\mathbb{R}_u^n)$.
\item[(b)] For every $L \in \mathcal{S}_u^p(\mathbb{S}^n)$, we have
\[\rho(g_u(L),v) = \tan\rho_u(L,v), \qquad v \in \mathbb{S}_u.  \]
\end{enumerate}
\end{lem}
{\it Proof.} Statement (a) is an immediate consequence of Lemma \ref{gnombasic1} (a) and (b). \linebreak
From Lemma \ref{gnombasic1} (a) and the definitions of radial and spherical radial functions, we obtain
    \begin{align*}
        \rho(g_u(L),v) &= \max\{\lambda \geq 0 : \lambda v \in g_u(L)\}\\
        &= \max\left\{\lambda \geq 0 :\frac{u+\lambda v}{\|u+\lambda v\|}
            = \frac{1}{\sqrt{1+\lambda^2}}\ u + \frac{\lambda}{\sqrt{1+\lambda^2}}\,v \in L\right\}\\
        &= \tan \max\{\alpha \in [0,\mbox{$\frac{\pi}{2}$}): u \cos \alpha + v \sin \alpha \in L\}\\
        &= \tan \rho_u(L,v)
    \end{align*}
which proves (b). \hfill $\blacksquare$

\pagebreak

By Lemma \ref{radrad} (b), a function $\rho: \S_u \to [0,\frac{\pi}{2})$ is the spherical radial function
of a star set $L \in \mathcal{S}_u^p(\mathbb{S}^n)$ if and only if the $-1$-homogeneous extension
of $\tan \rho$ to $\mathbb{R}_u^n$ is the radial function of a star set in $\mathbb{R}_u^n$.

\vspace{0.2cm}

We call a proper star set $L \in \mathcal{S}_u^p(\S^n)$ a \emph{(spherical) star body} with respect to $u \in \mathbb{S}^n$
if $\rho_u(L,\cdot)$ is positive and continuous. Clearly, every proper convex body $K \in \mathcal{K}_u^p(\mathbb{S}^n)$ containing $u$ in its interior
is a star body with respect to $u$. \linebreak In order to establish a counterpart to (\ref{polar17}), we recall that, for $K \in \mathcal{K}^p(\mathbb{S}^n)$ with non-empty interior,
the polar body $K^\circ \in \mathcal{K}^p(\mathbb{S}^n)$ is defined by
\[K^\circ = \{v \in \mathbb{S}^n: v \cdot w \leq 0 \mbox{ for all } w \in K\} = \S^n \backslash\, \mathrm{int}\, K_{\frac{\pi}{2}}.\]
Note that if $K \in \mathcal{K}_u^p(\mathbb{S}^n)$ contains $u$ in its interior, then $K^\circ \in \mathcal{K}_{-u}^p(\mathbb{S}^n)$ contains $-u$ in its interior.

\begin{prop} If $u \in \mathbb{S}^n$ and $K \in \mathcal{K}_u^p(\mathbb{S}^n)$ contains $u$ in its interior, then
\begin{equation} \label{gnompolar}
g_u(K)^* = g_{-u}(K^\circ)
\end{equation}
and
\begin{equation} \label{suppradpolar}
h_u(K,\cdot) + \rho_{-u}(K^\circ,\cdot) = \frac{\pi}{2}.
\end{equation}
\end{prop}
{\it Proof.} By the definitions of the Euclidean and spherical polar bodies and the gnomonic projection, we have
\begin{align*}
g_u(K)^* &= \{x \in \R^n_u: x \cdot y \leq 1 \mbox{ for all } y \in g_u(K)\} \\
         &= \{x \in \R^n_u: x \cdot g_u(w) \leq 1 \mbox{ for all } w \in K\} \\
         &= \{x \in \R^n_u: w \cdot x \leq w \cdot u \mbox{ for all } w \in K\}\\
         &= \left\{x \in \R^n_u: w \cdot \frac{x - u}{\|x - u\|} \leq 0 \mbox{ for all } w \in K\right\}\\
         &= g_{-u}\left(\{v \in \S^n: v \cdot w \leq 0 \mbox{ for all } w \in K\}\right) = g_{-u}(K^\circ).
\end{align*}
which proves (\ref{gnompolar}). Lemma \ref{suppsupp}, (\ref{polar17}), (\ref{gnompolar}), and Lemma \ref{radrad} (b), now yield
\begin{align*}
    \tan h_u(K,\cdot) = h(g_u(K),\cdot)
    = \frac{1}{\rho(g_u(K)^*,\cdot)}
    = \frac{1}{\rho(g_{-u}(K^\circ),\cdot)}
    = \frac{1}{\tan \rho_{-u}(K^\circ,\cdot)}
\end{align*}
which is equivalent to (\ref{suppradpolar}). \hfill $\blacksquare$

\vspace{0.3cm}

Using spherical radial functions, we define a metric $\widetilde{\gamma}_u$ on $\mathcal{S}_u^p(\mathbb{S}^n)$ by
\begin{align*}
\widetilde{\gamma}_u(K,L) = \sup_{v\in \S_u} |\rho_u(K,v)-\rho_u(L,v)|.
\end{align*}
Note that if $K, L \in \mathcal{S}_u^p(\mathbb{S}^n)$, then by (\ref{defraddist}) and Lemma \ref{radrad} (b),
\begin{align*}
\widetilde{\delta}(g_u(K),g_u(L)) = \sup_{v\in \S_u} |\tan \rho_u(K,v) - \tan \rho_u(L,v)|.
\end{align*}
Thus, from the continuity of the tangent we obtain the following.

\begin{theorem} \label{gnomhomstar}
The gnomonic projection is a homeomorphism between $(\mathcal{S}_u^p(\S^n),\widetilde{\gamma}_u)$ and $(\mathcal{S}(\R^n_u),\widetilde{\delta})$.
\end{theorem}

For fixed $u \in \mathbb{S}^n$ we call a binary operation $*: \mathcal{S}_u^p(\S^n) \times \mathcal{S}_u^p(\S^n) \rightarrow \mathcal{S}_u^p(\S^n)$
\emph{$u$-section covariant} if for all $k$-spheres $S$, $1 \leq k \leq n - 1$, with $u \in S$ and for all $K, L \in \mathcal{S}_u^p(\S^n)$, we have
\[(K \cap S) * (L \cap S) = (K * L) \cap S.  \]
The operation $*$ is called \emph{$u$-rotation covariant} if $(\vartheta K) * (\vartheta L) = \vartheta(K * L)$ 
for all $\vartheta \in \mathrm{SO}(n + 1)$ which fix $u$. Our next result is a version of Theorem \ref{thm:2} (or Theorem \ref{thm:2refined}, respectively)
in the setting of star sets.

\begin{theorem} \label{usecurot} For $u \in \mathbb{S}^n$, the gnomonic projection $g_u$ induces a one-to-one correspondence
between operations $*\!: \mathcal{S}_u^p(\S^n) \times \mathcal{S}_u^p(\S^n) \rightarrow \mathcal{S}_u^p(\S^n)$ which are \linebreak $u$-rotation and $u$-section covariant
and operations $\overline{*}\!: \mathcal{S}(\R^n_u) \times \mathcal{S}(\R^n_u) \rightarrow \mathcal{S}(\R^n_u)$ which are rotation and section covariant.
Moreover, any such operation $*$ is continuous if and only if $\overline{*}$ is continuous.
\end{theorem}
{\it Proof.} First assume that $*$ is $u$-rotation and $u$-section covariant and define an operation
$\overline{*}: \mathcal{S}(\R^n_u) \times \mathcal{S}(\R^n_u) \rightarrow \mathcal{S}(\R^n_u)$ by
\[K\,\overline{*}\,L = g_u(g_u^{-1}(K) * g_u^{-1}(L))   \]
for $K, L \in \mathcal{S}(\R_u^n)$. As in the proof of Theorem \ref{thm:2refined}, it follows that $\overline{*}$ is section covariant. The rotation covariance of $\overline{*}$ is
a consequence of the $u$-rotation covariance of $*$ and the fact that $\vartheta g_u(L) = g_u(\vartheta L)$ 
for all $L \in \mathcal{S}_u^p(\S^n)$ and $\vartheta \in \mathrm{SO}(n + 1)$ which fix $u$.

Conversely, if $\overline{*}: \mathcal{S}(\R^n_u) \times \mathcal{S}(\R^n_u) \rightarrow \mathcal{S}(\R^n_u)$ is rotation and section covariant, then
define $*: \mathcal{S}_u^p(\S^n) \times \mathcal{S}_u^p(\S^n) \rightarrow \mathcal{S}_u^p(\S^n)$ by
\[K * L = g_u^{-1}(g_u(K)\,\overline{*}\,g_u(L))  \]
for $K, L \in \mathcal{S}_u^p(\S^n)$. 

\pagebreak

As before, it is easy to show that $*$ is $u$-rotation and $u$-section covariant and,
by Theorem \ref{gnomhomstar}, the operation $*$ is continuous if and only if $\overline{*}$ is continuous.
\hfill $\blacksquare$

\vspace{0.3cm}

We conclude with a corollary to Theorem \ref{thm:starshaped} of Gardner, Hug, and Weil and Theorem \ref{usecurot}.

\begin{koro} For fixed $u \in \mathbb{S}^n$, an operation $*: \mathcal{S}_u^p(\S^n) \times \mathcal{S}_u^p(\S^n) \to \mathcal{S}_u^p(\S^n)$
is $u$-rotation and $u$-section covariant if and only if there exists a function $f: [0,\frac{\pi}{2})^4 \to [0,\frac{\pi}{2})$ such that,
for all $K, L \in \mathcal{S}_u^p(\S^n)$,
\begin{align*}
\rho_u(K*L,v) = f(\rho_u(K,-v),\rho_u(K,v),\rho_u(L,-v),\rho_u(L,v)), \qquad v \in \mathbb{S}_u.
\end{align*}
\end{koro}

\vspace{0.3cm}

\noindent {{\bf Acknowledgments} The work of the authors was
supported by the European Research Council (ERC) within the
project ``Isoperimetric Inequalities and Integral Geometry'',
Project number: 306445.

\begin{small}

Vienna University of Technology \par Institute of Discrete
Mathematics and Geometry \par Wiedner Hauptstra\ss e 8--10/1046
\par A--1040 Vienna, Austria

\vspace{0.2cm}

\par florian.besau@tuwien.ac.at \par franz.schuster@tuwien.ac.at

\end{small}

\end{document}